\documentclass[12pt]{article}
\usepackage{amssymb,amsmath,amsfonts,amsthm,amscd,latexsym,indentfirst}
\def\rb{\textrm{rb}}

\def\Span{\mathrm{Span}}

\def\Imm{\mathrm{Im}\,}

\begin{document}

\begin{center}
{\Large
Rota---Baxter operators of zero weight on simple Jordan algebra of Clifford type}

V. Gubarev
\end{center}

\begin{abstract}
It is proved that any Rota---Baxter operator of zero weight on
Jordan algebra of a nondegenerate bilinear symmetric form
is nilpotent of index less or equal three.
We state exact value of nilpotency index on simple
Jordan algebra of Clifford type over fields
$\mathbb{R}$, $\mathbb{C}$, and $\mathbb{Z}_p$.
For $\mathbb{Z}_p$, we essentially use the results
from number theory concerned quadratic residues and
Chevalley---Warning theorem.

\end{abstract}

\section{Introduction}

In 1933, the physicist P. Jordan introduced in quantum mechanics a notion
of Jordan algebra \cite{Jordan}. A Jordan algebra $J$ over a field $F$
of characteristic $\neq2$ is a commutative (nonassociative) algebra over
$F$ satisfying the Jordan identity $(x^2 y)x = x^2(yx)$.

Any associative algebra under the new product
$x\circ y = \frac{1}{2}(xy+yx)$
is a Jordan algebra. For example, a matrix algebra under $\circ$
turns to be a simple Jordan algebra. Moreover,
the set of self-adjoint matrices over $\mathbb{R}$, $\mathbb{C}$,
or $\mathbb{H}$ with the product $\circ$ is a Jordan algebra
(of Hermitian type).
The set of self-adjoint matrices of order three over the octonions
under $\circ$ is a 27-dimensional simple Jordan algebra (of Albert type).
Given a vector space $V$ with nondegenerate bilinear symmetric form $f$,
one can equip a space $F\oplus V$ with a structure of simple Jordan algebra
provided $\dim (V)\geq2$ (Clifford type).

The first classification theorem for Jordan algebras was obtained
by P. Jordan, J. von Neumann, and E. Wigner in 1934 \cite{JNW}.
In 1978--1983, E. Zelmanov was a leader of so called Russian revolution
\cite{McCrimmon} in theory of Jordan algebras.
In particular, he proved that any prime (including simple) Jordan algebra
without absolute zero divisors is either of Hermitian, or Albert, or
Clifford type \cite{Zelmanov3}.

Jordan algebras are deeply connected with Lie algebras and finite simple
groups (Kantor---K\"{o}cher---Tits construction,
The Freudenthal---Tits Magic Square) \cite{Jacobson,McCrimmon}.

Let us refer a reader to the excellent monographs on Jordan algebras
\cite{Jacobson,McCrimmon,Nearly}
written by N. Jacobson, K. McKrimmon, and K.A. Zhevlakov et al.

The notions of Rota---Baxter operator and Rota---Baxter
algebra were introduced by G. Baxter in 1960 \cite{Baxter60}
as a generalization of by part integration formula.
Since 1960s, J.-C. Rota, P. Cartier, L. Guo, C. Bai and others
studied combinatorial and algebraic properties of Rota---Baxter algebras
(see details in a unique monograph on the subject written by L. Guo
in 2012 \cite{GuoMonograph}).

Also, the deep connection of Rota---Baxter algebras with Yang---Baxter equation,
quantum field theory and Loday algebras was found
\cite{BelaDrin82,Semenov83,Aguiar00,QFT00,BBGN2012,Embedding}.

In 1935, C. Chevalley \cite{Chevalley} and E. Warning \cite{Warning}
proved two very close theorems which are now gathered
in one Chevalley---Warning theorem. This theorem states than over finite
field any system of polynomial equations
with sufficiently large number of variables has a solution.
From that time, a lot of connections and applications of
Chevalley---Warning theorem were found such that
combinatorial Nullstellensatz,
Erd\"{o}s---Ginzburg---Ziv theorem etc. (see, e.g. \cite{Clark}).

The goal of the current work is to investigate Rota---Baxter operators
of zero weight on simple Jordan algebra of a nondegenerate bilinear form.
At first, we prove that all such operators are nilpotent of index not greater than 3.
At second, we are interested on the exact value of nilpotency index
of Rota---Baxter operators of zero weight. It occurs that the nilpotency index
of Rota---Baxter operators of zero weight on simple Jordan algebra of Clifford type
depends on the dimension of algebra and the ground field. If the ground field
is finite, we essentially use Chevalley---Warning theorem and its corollaries.

We are always mentioning that characteristic of the ground field $F$ is not 2.

\section{Preliminaries}

\subsection{RB-operator}

Given an algebra $A$ and scalar $\lambda\in F$, where $F$ is a ground field,
a linear operator $R\colon A\rightarrow A$ is called
a Rota---Baxter operator (RB-operator, for short) on $A$ of the weight $\lambda$
if the following identity
\begin{equation}\label{RB}
R(x)R(y) = R( R(x)y + xR(y) + \lambda xy )
\end{equation}
holds for any $x,y\in A$.

Let us list some examples of RB-operators of zero weight (see, e.g., \cite{GuoMonograph}):

{\bf Example 1}.
Given an algebra $A$ of continuous
functions on $\mathbb{R}$, an integration operator
$R(f)(x) = \int\limits_{0}^x f(t)\,dt$
is an RB-operator on $A$ of zero weight.

{\bf Example 2}.
A linear map $R$ on polynomial algebra $F[x]$ defined as
$R(x^n) = \frac{x^{n+1}}{n+1}$ is an RB-operator on $F[x]$ of zero weight.

{\bf Example 3}.
Given an invertible derivation $d$ on an algebra $A$,
$d^{-1}$ is an RB-operator on $A$ of zero weight.

{\bf Lemma 1}.
Let $A$ be an unital algebra, $P$ be an RB-operator on $A$ of zero weight.

a) $1\not\in \Imm P$;

b) if $A$ is simple finite-dimensional algebra and $\dim A>1$, then $\dim \ker P\geq2$;

c) if $P(1)\in F$, then $P(1) = 0$, $P^2 = 0$, and $\Imm P\subset \ker P$;

d) $P(1)P(1) = 2P^2(1)$.

{\sc Proof}.
a) Assume there exists $x$ such that $P(x) = 1$. We have
$$
1 = 1\cdot1 = P(x)P(x)
 = P(P(x)x + xP(x)) = 2P(x) = 2\cdot1,
$$
a contradiction.

b) Let $A$ be simple, $\dim A = n < \infty$.
Suppose that $\dim \Imm P = n-1$. Given a nonzero $x\in \ker P$,
we compute
$$
0 = P(x)P(y) = P(xP(y)),\quad
0 = P(y)P(x) = P(P(y)x).
$$
Therefore, $\ker P$ is closed under multiplication
on $\Imm P$. As $A = \Span\{1,\Imm P\}$,
we have a proper ideal $\ker P$ in $A$,
this contradicts simplicity of $A$.

c) From a), $P(1) = 0$. All other follows from
$$
0 = P(1)P(x) = P(P(1)x + 1\cdot P(x)) = P(P(x)).
$$

d) Follows from \eqref{RB} directly.

{\bf Lemma 2}.
Let an algebra $A$ equal $B \oplus C$ (as vector spaces), where $B$ is a subalgebra
of $A$ with an RB-operator $P$ of zero weight,
$BC,CB\subseteq \ker P\oplus C$.
Then the linear operator $R$ on $A$ defined as
$R(b+c) = P(b)$, $b\in B$, $c\in C$,
is an RB-operator of weight~0.

{\sc Proof}. Straightforward.

\subsection{Chevalley---Warning theorem}

For our purposes, let us write down the following version
of Chevalley---Warning theorem:

{\bf Chevalley---Warning theorem}.
Let $F$ be a finite field and
$f_i\in F[x_1,\ldots,x_n]$, $i=1,\ldots,r$, be  homogeneous
polynomials  of  degree $d_i$ respectively.
If $d_1+\ldots+d_r<n$, then  there  is a nonzero solution to
$f_1(x) = \ldots = f_r(x) = 0$.

We also need the following result

{\bf Proposition} \cite{Clark}.
Let $F$ be a field of characteristic different from 2 in which each
quadratic form in three variables has a nontrivial solution.
Then, for any $a,b,c\in F^*$, there exist
$x,y\in F$ such that $ax^2+by^2=c$.

For a finite field, Proposition could be derived from \cite{Heath-Brown}
as a corollary of Chevalley---Warning theorem.

\subsection{Simple Jordan algebra of Clifford type}

Let us introduce simple $(n+1)$-dimensional Jordan algebra
$J_{n+1}(f) = F\cdot1\oplus V$
of a nondegenerate bilinear symmetric form $f$ acting
on $n$-dimensional vector space $V$, $n\geq2$.
The product on $J_{n+1}(f)$ is the following:
$$
(\alpha\cdot1+v)(\beta\cdot1 +u)
 = (\alpha\beta+f(v,u))\cdot1 + (\alpha u + \beta v).
$$
Let us choice such base $e_1$, $e_2$, \ldots, $e_n$ of $V$ that
the matrix of the form $f$ in this base is diagonal with
elements $d_1,d_2,\ldots,d_n$ on the main diagonal.
As $f$ is nondegenerate, $d_i\neq0$ for each $i$.
Let us identify the form $f$ with the elements $d_i$:
$f = (d_1,\ldots,d_n)$.
By $(v,u)$ for $v = \sum\limits_{i=1}^n v_i e_i$,
$u = \sum\limits_{i=1}^n u_i e_i$,
we denote the sum $\sum\limits_{i=1}^n d_i v_i u_i$.

It is well-known \cite{Nearly} that $J_{n+1}(f)$
is a quadratic algebra, it means that any element
$x\in J_{n+1}(f)$ satisfies a quadratic equation
$x^2 - t(x)x + n(x)1 = 0$ for some $t(x),n(x)\in F$.

\section{Bound theorem and 3-dimensional case}

{\bf Theorem 1}.
Let $J$ be a (not necessary simple) Jordan algebra of a bilinear form,
then $\rb(J)\leq 3$.

{\sc Proof}.
Let us consider $R(1) \neq 0$, else $R^2 = 0$ by Lemma 1,\,c).
By Lemma 1,\,a) we have $n(R(x)) = 0$ for each $x\in A$.
As $R(1)R(1) = tR(1)$ for $t = t(R(1))\in F$, there are two cases:
$t = 0$ or $t\neq0$.

Suppose that $t = 0$. Let us compute by Lemma 1,\,d) for $x\in V$
\begin{equation}\label{Jord0}
0 = R^2(1)R(x)
  = R( R^2(1)x + R(1)R(x) )
  = R^2( R(1)x + R(x) )
  = R^2(R(1)x) + R^3(x).
\end{equation}
As $R(1)\in V$, we have that $R(1)x\in F$ for all $x\in V$.
So $R^2(R(1)x) = 0 = R^3(x)$, $x\in V$.
Together with $R^3(1) = R(R^2(1)) = 0$
we obtain that $R^3 = 0$ on the whole $A$.

Let $t\neq0$ and $R(1) = \alpha\cdot1+a$ for $\alpha=t/2\neq0$, $a\in V$.
From $n(R(1)) = 0$ we get $a^2 = \alpha^2\cdot1$. At first,
$$
2\alpha R(1)
 = R(1)R(1)
 = 2R^2(1) = 2R(\alpha\cdot1+a)
 = 2\alpha R(1) + 2R(a),
$$
so $R(a) = 0$.
At second,
$$
0 = R(1)R(a)
  = R( R(1)a + R(a) )
  = R( R(1)a)
  = R( (\alpha\cdot1+a)a)
  = \alpha^2 R(1),
$$
a contradiction.

In light of Theorem 1, define nilpotency index $\rb(J_{n+1}(f))$
as a minimal natural number $s$ such that any RB-operator of
zero weight on $J_{n+1}(f)$ is inlpotent in the power~$s$.

{\bf Statement}.
All RB-operators of zero weight on simple Jordan algebra
$J_3(f)$, $f = (d_1,d_2)$, are the following
\begin{equation}\label{J3-0}
R(1) = 0,\quad R(e_1) = k(a\cdot1+be_1+ce_2),\quad
R(e_2) = l(a\cdot1+be_1+ce_2)
\end{equation}
for such $a,b,c,k,l\in F$ that $a^2-d_1 b^2-d_2 c^2 =0$ and $kb+lc = 0$.
Moreover, $R^2 = 0$.

{\sc Proof}.
Let $R$ be a nonzero RB-operator on $A$ of zero weight.
By Lemma 1,\,b), $\dim\ker R = 2$ and $\dim\Imm R = 1$.
Proving $R(1) = 0$, we have exactly \eqref{J3-0}.
The condition $a^2-d_1 b^2-d_2 c^2 = 0$ follows from the fact
that $\Imm R$ is a subalgebra of $A$ and $kb + lc = 0$
follows from \eqref{RB}.

Suppose that $R(1)\neq0$. By Theorem 1,
$R^2 = 0$, $\Imm R\subset \ker R$, and $R(1) = ae_1+be_2$ for some $a,b\in F$.
Let $x = \chi1+ke_1+le_2\in \ker R$ be such element that
the vectors $R(1)$ and $ke_1+le_2$ are linear independent. From
$$
0 = R(1)R(x)
  = R(R(1)x) = R(\chi R(1)+c\cdot1) = cR(1),
$$
where $c = d_1ak+d_2bl$, we get $c = 0$.
Together with $d_1a^2 + d_2b^2 = 0$, a condition
of $\Imm R$ being a subalgebra of $A$, we have $al = bk$,
i.e., $R(1)$ and $ke_1+le_2$ are linear dependent,
a contradiction.

{\bf Corollary 1}.
Let $F$ be an algebraically closed or finite field,
then $\rb(J_3) = 2$.

{\sc Proof}.
The main question is to find a nonzero solution $(a,b,c)$
of the equation $a^2 - d_1b^2 - d_2 c^2 = 0$.
Over algebraically closed field it is enough to consider
$c = 1$, $b = 0$ and $a = \pm\sqrt{d_2}$. For a finite field,
we apply Chevalley---Warning theorem.

The condition $kb + lc = 0$ could be easily satisfied.
At least one of the numbers $b,c$ is nonzero.
Suppose that $b\neq0$, choose any $l\neq0$ and put $k = -lc/b$.

{\bf Remark 1}.
The statement of Corollary 1 could be not true for other fields.
For example, in the case $F = \mathbb{R}$ and $f = (1,1)$,
we have $\rb(J_3(f)) = 1$.

\section{General case}

{\bf Lemma 3}.
Let $J_{n+1}(f)$ be $(n+1)$-dimensional simple Jordan algebra
of the form $f = (d_1,\ldots,d_n)$.
Let $R(1) = k_1 e_1 + \ldots + k_n e_n$,
$R(e_i) = \alpha_0\cdot1 + \alpha_1 e_1 + \ldots + \alpha_n e_n$,
$R(e_j) = \beta_0\cdot1 + \beta_1 e_1 + \ldots + \beta_n e_n$,
$k_a,\alpha_b,\beta_c\in F$. Suppose, $R(1)\neq0$. Then

a) $R(1)R(1) = R^2(1) = (R(1),R(1)) = 0$;

b) $R(1)R(e_i) = \alpha_0 R(1)$ and
$(R(1),R(e_i)-\alpha_0\cdot1) = 0$;

c) $R^2(e_i) = (\alpha_0-d_i k_i)R(1)$;

d) $R(e_i)R(e_j) = \alpha_0 R(e_j) + \beta_0 R(e_i)$,
$(R(e_i)-\alpha_0\cdot1,R(e_j)-\beta_0\cdot1) = \alpha_0\beta_0$,
and $d_j\alpha_j + d_i\beta_i = 0$.

e) $R(e_i)R(e_i) = 2\alpha_0 R(e_i)$,
$(R(e_i)-\alpha_0\cdot1,R(e_i)-\alpha_0\cdot1) = \alpha^2_0$,
and $\alpha_i = 0$.

{\sc Proof}.
a) Follows from the proof of Theorem 1 and Lemma 1,\,d).

b) By Lemma 1,\,b),
$R(1)R(e_i) = \alpha_0 R(1) + (R(1),R(e_i)-\alpha_0\cdot1)\cdot1
= \alpha_0 R(1)$ and $(R(1),R(e_i)-\alpha_0\cdot1) = 0$.

c) Follows from b) and application of \eqref{RB} for $R(1)R(e_i)$.

d) Let $v = R(e_i)-\alpha_0\cdot1$, $u = R(e_j)-\beta_0\cdot1$.
From one hand,
\begin{equation}\label{RB:Jordan:productL}
R(e_i)R(e_j)
 = (\alpha_0\beta_0+(v,u))\cdot1
 + \alpha_0 u + \beta_0 v
 = (-\alpha_0\beta_0+(v,u))\cdot1
 + \alpha_0 R(e_j) + \beta_0 R(e_i).
\end{equation}
From another hand,
\begin{multline}\label{RB:Jordan:productR}
R(R(e_i)e_j+e_iR(e_j))
 = R(\alpha_0 e_j + d_j\alpha_j\cdot1) +
  R(\beta_0 e_i + d_i\beta_i\cdot1) \\
  = (d_j\alpha_j + d_i\beta_i)R(1) +
 \alpha_0 R(e_j) + \beta_0 R(e_i).
\end{multline}
Comparing \eqref{RB:Jordan:productL} and \eqref{RB:Jordan:productR},
we have done.

e) Follows from d) for $i = j$.

{\bf Remark 2}.
Suppose that $R$ is a such RB-operator on $J_{n+1}(f)$ that $R^2\neq0$.
Moreover, $\sqrt{d_j}\in F$. Define $e_i' = e_i/\sqrt{d_i}$.
Let $A$ be a matrix of $R$ in the base $1,e_1',\ldots,e_n'$
and $M$ its submatrix of the size $n\times n$ formed by all
rows and columns except the first ones. By Lemma 3,\,d), we have
that $M$ is skew-symmetric.

{\bf Example 4}.
Let $J_{n+1}(f)$ be $(n+1)$-dimensional simple Jordan algebra
of the form $f = (d_1,\ldots,d_n)$.
Let $l_1,\ldots,l_p,k_{p+1},\ldots,k_n\in F$ be such that
$k_j\neq0$ for some $j$ and
\begin{equation}\label{RB:JordanExmpConst}
\sum\limits_{i=1}^p d_i l_i^2 = 1,\quad
\sum\limits_{j=p+1}^n d_j k_j^2 = 0.
\end{equation}
Then the the following linear map on $J_{n+1}(f)$
\begin{equation}\label{RB:JordanExmpFormA}
R(e_i) = \begin{cases}
k_{p+1} e_{p+1}+\ldots + k_n e_n, & i = 0, \\
d_i l_i(k_{p+1} e_{p+1}+\ldots + k_n e_n), & 1\leq i\leq p, \\
-d_i k_i(1+l_1 e_1+\ldots+l_p e_p), & p+1\leq i\leq n,
\end{cases}
\end{equation}
is an RB-operator of zero weight.
Moreover, $R^2\neq 0$ and $\rb(J_{n+1}(f)) = 3$.

{\bf Example 5}.
Let $J_4(f)$ be $4$-dimensional simple Jordan algebra
of the form $f = (d_1,d_2,d_3)$.
Moreover, there exist a solution $x_0$ of the equation
$x^2 + d_1 d_2 d_3 = 0$ and such $k_1,k_2,k_3\in F^*$ that
$d_1 k_1^2 + d_2 k_2^2 + d_3 k_3^2 = 0$.
Then the the following linear map on $J_4(f)$
\begin{equation}\label{RB:JordanExmpFormB}
\begin{gathered}
R(1) = d_1 k_1^2 e_1' + d_2 k_2^2 e_2' + d_3 k_3^2 e_3', \quad
R(e_1') = -1 + \frac{\lambda}{d_1 k_1^2}(e_2'-e_3'), \\
R(e_2') = -1 + \frac{\lambda}{d_2 k_2^2}(-e_1'+e_3'), \quad
R(e_1') = -1 + \frac{\lambda}{d_3 k_3^2}(e_1'-e_2'),
\end{gathered}
\end{equation}
where $\lambda = k_1 k_2 k_3 x_0$, $e' = e_i/(k_i d_i)$,
is an RB-operator of zero weight and $\rb(J_4(f)) = 3$.

{\bf Theorem 2}.
Let $F$ be a field, $J_{n+1}(f)$ is a simple Jordan algebra
over $F$ of the form $f = (d_1,\ldots,d_n)$, $n\geq3$.

a) If $i\in F$ and there exist $j_1,j_2,j_3\in\{1,\ldots,n\}$ such that
$\sqrt{d_{j_s}}\in F$, $s=1,2,3$, then $\rb(J_{n+1}(f)) = 3$.

b) For $F = \mathbb{R}$, we have
$$
\rb(J_{n+1}(f)) = \begin{cases}
1, & d_1,\ldots,d_n<0, \\
2, & d_1,\ldots,d_n>0\ \mbox{or}\ d_1,\ldots,d_{i-1},d_{i+1},\ldots,d_n<0,\,d_i>0, \\
3, & else.
\end{cases}
$$

c) If $F$ is a finite field, then $\rb(J_k(f)) = 3$ for all $k\geq 6$.

d) If $F = \mathbb{Z}_p$ for a prime $p$ of the form $4t+1$, then we have
$$
\rb(J_k(f)) = \begin{cases}
2, & k=3\ \mbox{or }k=4,\,\mbox{even number of residues through }d_1,d_2,d_3, \\
3, & k\geq5\ \mbox{or }k=4,\,\mbox{odd number of residues through }d_1,d_2,d_3.
\end{cases}
$$

e) If $F = \mathbb{Z}_p$ for a prime $p$ of the form $4t+3$, then we have
$$
\rb(J_k(f)) = \begin{cases}
2, & k=3\ \mbox{or }k=4,\,\mbox{odd number of residues through }d_1,d_2,d_3, \\
3, & k\geq5\ \mbox{or }k=4,\,\mbox{even number of residues through }d_1,d_2,d_3.
\end{cases}
$$

{\sc Proof}.
a) Let $n = 3$, by Example 4 the following linear map $P$ on $J_4(f)$
\begin{equation}\label{RB:Jordan-big-C}
\begin{gathered}
P(1) = \frac{e_2}{\sqrt{d_2}} + \frac{ie_3}{\sqrt{d_3}}, \quad
P(e_1) = \sqrt{d_1}\left(\frac{e_2}{\sqrt{d_2}} + \frac{ie_3}{\sqrt{d_3}}\right), \\
P(e_2) = -\sqrt{d_2}\left(1+\frac{e_1}{\sqrt{d_1}}\right), \quad
P(e_3) = -\sqrt{d_3}\left(1+\frac{e_1}{\sqrt{d_1}}\right)
\end{gathered}
\end{equation}
is an RB-operator of zero weight and $P^2\neq0$.

Let $n > 3$ and suppose that $j_s = s$, $s=1,2,3$.
Define a linear map $R$ on $J$ by Lemma 2
for $B = \Span\{1,e_1,e_2,e_3\}$ with $P$ and $C = \Span\{e_4,\ldots,e_n\}$.

b) If all $d_i$ are negative, then $R = 0$ by Lemma 3,\,a),\,d).
If all $d_i$ are positive, then $R(1) = 0$ by Lemma 3,\,a).
Thus, $\rb(J_{n+1}(f)) = 2$ by Statement.
Suppose that $d_2,\ldots,d_n<0$ and $d_1>0$.
Let $R$ be a such RB-operator on $J_{n+1}(f)$ that $R(1)\neq0$.
Denote $R(e_i) = \sum\limits_{j=0}^n\alpha_{ij}$, $i=1,\ldots,n$,
where $e_0 = 1$. By Lemma 3,\,e), $\alpha_{11} = 0$ and $R(e_1) = 0$.
From Lemma 3,\,d), we get $\alpha_{i1} = 0$ for all $i=1,\ldots,n$.
So, $R(e_i) = 0$, $i=1,\ldots,n$. Hence, $R^2 = 0$.
Thus, $\rb(J_{n+1}(f)) = 2$ by Statement.

Finally, let $d_1,d_2>0$, $d_3<0$, all other $d_i$ are of some sign.
We apply Example 4 for $p = 1$ to obtain $\rb(J_{n+1}(f)) = 3$.

c) Let $n = 5$ and $p = 2$.
By Chevalley---Warning theorem, the second equation \eqref{RB:JordanExmpConst}
has a nonzero solution $(k_3,k_4,k_5)$.
By Proposition, the first equation of \eqref{RB:JordanExmpConst} also has a solution.
Applying Example 4, we are done.

We deal with $n > 5$ in analogous way as in a).

d) Let $n = 3$. Suppose, $d_1$ is a quadratic residue modulo $p$
and $d_2,d_3$ are simultaneously quadratic residues or nonresidues.
By Example 4, for $p = 1$, we have $\rb(J_4(f)) = 3$.

Suppose, there is an odd number of residues through $d_1,d_2,d_3$.
To the contrary, let $R$ be an RB-operator of zero weight on $J_4(f)$
such that $R^2\neq 0$. By Lemma 1,\,c), $R(1) = k_1e_1+k_2e_2+k_3e_3\neq 0$.
Without loss of generality we assume that $k_3\neq0$.
By Lemma 1,\,b) and assumption that $R^2\neq 0$, we have
$\dim\Imm R = \dim\ker R = 2$.

Consider $R(e_1) = \alpha_0\cdot1+\alpha_1e_1+\alpha_2e_2+\alpha_3e_3$.
By Lemma 3,\,e), $\alpha_1 = 0$.
By Lemma 3,\,b), $d_2k_2\alpha_2+d_3k_3\alpha_3 = 0$ and
by Lemma 3,\,e), $\alpha_0^2 = d_2\alpha_2^2+d_3\alpha_3^2$.
From the last two equalities and Lemma 3,\,a), we conclude that
\begin{multline}\label{RB:Jordan:Kontra}
\alpha_0^2
 = d_2\alpha_2^2\left(1+\frac{d_3}{d_2}\left(\frac{d_2k_2}{d_3k_3}\right)^2\right)
 = d_2\alpha_2^2\left(1+\frac{d_2k_2^2}{d_3k_3^2}\right) \\
 = \frac{d_2\alpha_2^2}{d_3k_3^2}(d_2k_2^2+d_3k_3^2)
 = -\frac{d_2d_1k_1^2\alpha_2^2}{d_3k_3^2}.
\end{multline}

If all $k_i$ are nonzero, then by \eqref{RB:Jordan:Kontra},
$(d_3k_3\alpha_0,k_1\alpha_2)$ is a solution of the equation
$x^2+d_1d_2d_3y^2 = 0$, it is a contradiction with a choice
of $d_1,d_2,d_3$.

Let $k_1 = 0$ and $k_2,k_3$ be nonzero.
Suppose that $R(1)$ and $R(e_1)$ form a base of $\Imm R$.
By \eqref{RB:Jordan:Kontra}, $\alpha_0 = 0$ and
$\Imm R = \Span\{e_2,e_3\}$. So, coordinates of $R(e_2)$ and $R(e_3)$
on $1,e_1$ are zero. Applying Lemma 3,\,d), we have $R(e_1) = 0$,
a contradiction.

Suppose that $R(1)$ and $R(e_2) = \beta_0\cdot1+\beta_1e_1+\beta_2e_2+\beta_3e_3$
form a base of $\Imm R$. By Lemma 3,\,e), $\beta_2 = 0$.
As above, we have $\beta_3 = 0$. So, the equations
$x^2-d_1y^2 = 0$ and $d_2z^2 + d_3t^2 = 0$
have nonzero solutions $(\beta_0,\beta_1)$ and $(k_2,k_3)$
respectively, it contradicts a choice of $d_1,d_2,d_3$.

Let $n\geq 4$, then by the properties
of quadratic residues we have that at least one of the equations
$d_i x^2 + d_j y^2 = 0$ for $i,j\in\{1,2,3\}$, $i\neq j$,
has a nonzero solution. The rest follows due to the arguments from b).

e) Let $n = 3$. If $d_1,d_2$ are quadratic residues modulo $p$
and $d_3$ is not, we have done by Example 4 for $p = 1$.
If $d_1,d_2,d_3$ are quadratic nonresidues modulo $p$,
we have $\rb(J_4(f)) = 3$ by Example 5.
Indeed, by Chevalley---Warning theorem, the second equation
\eqref{RB:JordanExmpConst} has a nonzero solution $(k_1,k_2,k_3)$.
The equation $x^2 + d_1 d_2 d_3 = 0$ has a solution by
the properties of residues. The proof in the case
when $R(1)$, $R(e_3)$ form a base of $\Imm R$ is analogous.

Supposing that there is an even number of residues through $d_1,d_2,d_3$,
we have
\linebreak
$\rb(J_4(f)) = 2$ by analogous arguments as in d).

For $n\geq4$, if not all $d_i$ are quadratic residues modulo $p$,
we have $\rb(J_{n+1}(f)) = 3$ by extending by Lemma 2
an RB-operator from $n = 3$. When all $d_i$ are quadratic residues modulo $p$,
we have done by Example 4 for $p = 1$.

{\bf Corollary 2}.
Let $F$ be an algebraically closed field and $J_{n+1}(f)$
be a simple Jordan algebra of a bilinear form $f$ over $F$.
We have
$\rb(J_k(f)) = \begin{cases}
2, & k=3, \\
3, & k\geq4.
\end{cases}$

{\bf Example 6}.
Let $F = \mathbb{Z}_7$, $J_4(f)$ be a $4$-dimensional simple Jordan algebra
of the form $f = (-1,-1,-1)$. The following RB-operator $R$ on $J_4(f)$
is defined by Example 5:
$$
R(1) = e_1 + 2e_2 + 3e_3, \quad
R(e_1) = 1 + 4e_2 + 2e_3, \quad
R(e_2) = 2 + 3e_1 + 6e_3, \quad
R(e_3) = 3 + 5e_1 + e_2.
$$

For investigating $\rb(J_k(f))$ over the field $\mathbb{Q}$,
Hasse---Minkowski theorem could be used. It says that a quadratic form $S$
has a rational nonzero solution if and only $S$ has a nonzero solution
over the real numbers and the $p$-adic numbers for every prime $p$.

{\bf Example 7}.
Let $F = \mathbb{Q}$, $J_4(f)$ be a $4$-dimensional simple Jordan algebra
of the form $f = (1,-3,1)$. Then $\rb(J_4(f)) = 2$.
Indeed, by Statement, $\rb(J_4(f)) \geq 2$.
Suppose, there exists an RB-operator $R$ on $J_4(f)$
such that $R^2\neq0$. Thus, $R(1) = k_1 e_1 + k_2 e_2 + k_3 e_3 \neq0$
and $k_1^2 - 3k_2^2 + k_3^2 = 0$. This equation hasn't nonzero
solution because the Pell equation $x^2 - 3y^2 = -1$ has no integer solutions.

\noindent
Gubarev Vsevolod \\
Sobolev Institute of Mathematics of the SB RAS \\
Acad. Koptyug ave., 4 \\
Novosibirsk State University \\
Pirogova str., 2 \\
Novosibirsk, Russia, 630090\\
{\it E-mail: wsewolod89@gmail.com}


\begin{thebibliography}{99}
\bibitem{Aguiar00}
Aguiar M.
Pre-Poisson algebras // Lett. Math. Phys. 2000. Vol. 54. P. 263--277.

\bibitem{BBGN2012}
Bai C., Bellier O., Guo L., Ni X.
Splitting of operations, Manin products, and Rota---Baxter operators //
Int. Math. Res. Notes. 2013. Vol. 3. P. 485--524.

\bibitem{Baxter60}
Baxter G.
An analytic problem whose solution follows from a simple
algebraic identity // Pacific J. Math. 1960. Vol. 10. P. 731--742.

\bibitem{BelaDrin82}
Belavin A.A., Drinfel'd V.G.
Solutions of the classical Yang-Baxter equation for
simple Lie algebras // Funct. Anal. Appl. 1982. Vol. 16, N 3. P. 159--180.

\bibitem{Chevalley}
Chevalley C.
D\'{e}monstration d'une hypoth\`{e}se de M. Artin //
Abh. Math. Sem. Univ. Hamburg. 1935. Vol. 11. P. 73--75.

\bibitem{Clark}
Clark P.L.
The Chevalley-Warning Theorem (Featuring... the Erd\"{o}s-Ginzburg-Ziv Theorem),
http://math.uga.edu/~pete/4400ChevalleyWarning.pdf

\bibitem{QFT00}
Connes A., Kreimer D.
Renormalization in quantum field theory and the Rie\-mann-Hilbert problem. I.
The Hopf algebra structure of graphs and the main theorem //
Comm. Math. Phys. 2000. Vol. 210, No. 1. P. 249--273.

\bibitem{Embedding}
Gubarev V., Kolesnikov P.
Embedding of dendriform algebras into Rota---Baxter algebras //
Cent. Eur. J. Math. 2013. Vol. 11, No. 2. P. 226--245.

\bibitem{GuoMonograph}
Guo L.
An Introduction to Rota---Baxter Algebra.
Surveys of Modern Mathematics. Vol. 4.
Somerville, MA: International Press; Beijing: Higher education press, 2012.

\bibitem{Heath-Brown}
Heath-Brown D.R.
A Note on the Chevalley---Warning Theorems //
Russian Math. Surveys. 2011. Vol. 66, no. 2. P. 427--436.

\bibitem{Jacobson}
Jacobson N.
Structure and representations of Jordan algebras.
Providence, USA: AMS Colloquium Publications, vol. XXXIX, 1968.

\bibitem{Jordan}
Jordan P.
\"{U}ber Verallgemeinerungsmoglichkeiten des Formalismus der Quantenmechanik //
Nachr. Akad. Wiss. Gottingen. Math. Phys. Kl. I, vol. 41. 1933. P. 209--217.

\bibitem{JNW}
Jordan P., von Neumann J., Wigner E.
On an Algebraic Generalization of the Quantum Mechanical Formalism //
Ann. Math. 1934. Vol. 35, iss. 1. P. 29--64.

\bibitem{McCrimmon}
McCrimmon K.
A taste of Jordan algebras. N.-Y.: Springer-Verlag New York Inc., 2004.

\bibitem{Semenov83}
Semenov-Tyan-Shanskii M.A.
What is a classical r-matrix? // Funct. Anal. Appl. 1983. Vol. 17, N 4. P. 259--272.

\bibitem{Warning}
Warning E.
Bemerkung zur vorstehenden Arbeit von Herrn Chevalley //
Abh. Math. Sem. Hamburg. 1935. Vol. 11. P. 76--83.

\bibitem{Zelmanov3}
Zelmanov E.
On prime Jordan algebras II // Siberian Math. J. 1983. Vol. 24. P. 73--85.

\bibitem{Nearly}
Zhevlakov K.A., Slin'ko A.M., Shestakov I.P., and Shirshov A.I.
Rings that are nearly associative. N.-Y.: Academic Press, 1982.

\end{thebibliography}
\end{document}